\documentclass[a4paper,14pt]{article}

\usepackage{multicol}
\usepackage{amsthm}

\usepackage[dvips]{graphicx,color}
\usepackage{amssymb}
\usepackage{amsmath}

\begin{document}

\fontsize{14pt}{16.5pt}\selectfont

\begin{center}
\bf{A not metrizable space filled with any $n$ mutually disjoint self-similar spaces}
\end{center}

\fontsize{12pt}{11pt}\selectfont
\begin{center}
Akihiko Kitada$^{1}$, Shousuke Ohmori$^{2}$, Tomoyuki Yamamoto$^{1,2,*}$\\ 
\end{center}

\noindent
$^1$\it{Institute of Condensed-Matter Science, Comprehensive Resaerch Organization, Waseda University,
3-4-1 Okubo, Shinjuku-ku, Tokyo 169-8555, Japan}\\
$^2$\it{Faculty of Science and Engineering, Waseda~University, 3-4-1 Okubo, Shinjuku-ku, Tokyo 169-8555, Japan}\\
*corresponding author: tymmt@waseda.jp\\
~~\\
\rm
\fontsize{10pt}{11pt}\selectfont\noindent

\noindent
{\bf Abstract}\\
A $\Lambda (Card\Lambda \succ \aleph)$-product space of $\{0,1\}$ has a partition $\{X_1,\dots,X_n\}$ for any $n$ a decomposition space of each $X_i$ of which is self-similar.(February 16, 2015)\\


\section{Introduction}

Let $(X, \tau)=(\{0,1\}^\Lambda ,\tau_0^\Lambda ),~Card \Lambda \succ \aleph$
(aleph) be the $\Lambda -$product space of $(\{0,1\},\tau_0)$ where $\tau_0$ is a 
discrete topology for $\{0,1\}$. The topological space $(X,\tau)$ need not to be metrizable, 
then, the relation $Card \Lambda \succ \aleph $ (aleph) holds and since the space $X$ is not metrizable, $X$ is no more the so-called Cantor set. In the present report, we 
will mathematically confirm the existence of a partition $\{X_1,\dots,X_n\}$ for any $n$ of $X$ ({\bf 2-a)}), $X_i \in 
{(\tau\cap \Im)}$\footnote{The notaion $\tau\cap \Im$ denotes the set of all closed and open sets (clopen sets) of $(X,\tau)$.}$-\{\phi\}$, a decomposition space (i.e, a space of equivalence 
classes) of each element $X_i=(X_i,\tau_{X_i})$\footnote{$\tau_A=\{u\cap A;u\in \tau\}$ for $A\subset X$.} of which is self-similar. 
Namely, all points in $X_i$ are classified into equivalence classes and then the equivalence classes form a self-similar structure. Each equivalence class is a point of the decomposition space of $X_i$.

Some easily verified and well known statements are summarized in the next section for the preliminaries.

\section{Preliminaries}
\begin{description}
\item[2-a)] A family $\{X_1,\dots,X_n\}$ of subspace $X_i$ of $X$ is a partition of $X$ provided that $X_i\cap X_j=\phi, i\not = j$ and $X_1\cup \cdot \cdot \cdot \cup X_n=X.$
\item[2-b)] Any zero-dimensional (0-dim), perfect, T$_0$(necessarily T$_2$)-space $(X,\tau)$ has a partition $\{X_1,\dots,X_n\}, X_i\in(\tau\cap \Im)-\{\phi\}$ for any $n$. Each subspace $(X_i,\tau_{X_i})$ is a 0-dim, perfect, T$_2$-space.
\item[2-c)] Let $X$ be a 0-dim, perfect, compact T$_2$-space. Then, for {\it any} compact metric space $Y$, there exists a continuous map $f$ from $X$ onto $Y$.
\item[2-d)] If $f:(X,\tau)\rightarrow (Y,\tau')$ is a quotient map, then the map $h:(Y,\tau')\rightarrow (\mathcal{D}_f,\tau(\mathcal{D}_f)), y\mapsto f^{-1}(y)$ is a homeomorphism. Here, the decomposition $\mathcal{D}_f$ of $X$ and the decomposition topology $\tau(\mathcal{D}_f)$ are given by $\mathcal{D}_f=\{f^{-1}(y)\subset X;y\in Y\}$ and $\tau(\mathcal{D}_f)=\{\mathcal{U}\subset \mathcal{D}_f;\bigcup \mathcal{U}\in\tau \}$, respectively. $\mathcal{D}_f$ forms a partition of $X$.
\item[2-e)] A space which is homeomorphic to a self-similar space is self-similar.
\item[2-f)] Let $(X,\tau)$ be a 0-dim, perfect, compact, not-metrizable T$_2$-space. From the Urysohn's metrization theorem $(X,\tau)$ is not second countable. In the partition $\{(X_1,\tau_{X_1}),\dots, (X_n,\tau_{X_n})\}$ of $(X,\tau)$, there exists a number $i_0$ such that $(X_{i_0},\tau_{X_{i_0}})$ is not second countable. The subspace $(X_{i_0},\tau_{X_{i_0}})$ is a 0-dim, perfect, compact not-metrizable T$_2$-space.
\end{description}

\section{A partition $\{X_1,\dots,X_n\}$ of $(\{0,1\}^\Lambda ,\tau_0^\Lambda )$}

Since $(X,\tau)=(\{0,1\}^\Lambda ,\tau_0^\Lambda ),\tau_0=2^{\{0,1\}},Card 
\Lambda \succ \aleph$ is easily verified to be a 0-dim, perfect, compact T$_2$-space, from 
{\bf 2-b)} there exists a partition $\{(X_1,\tau_{X_1}),\dots,(X_n,\tau_{X_n})\}$ of $(X,\tau)$ where each subspace $(X_i,\tau_{X_i})$ is a 0-dim, perfect, 
compact T$_2$-space.\footnote{As mentioned in {\bf 
2-f)}, there exists 0-dim, perfect, compact, not metrizable T$_2$-subspace $(X_{i_0},\tau_{X_{i_0}})$.} From {\bf 2-c)}, there exists a continuous map $f_i$ 
from $(X_i,\tau_{X_i})$ onto any compact, self-similar metric space $(Y,\tau_d)$. Since $(X_i,\tau_{X_i})$ is a compact space and $Y$ is a T$_2$-space, the map $f_i:(X_i,\tau_{X_i})\rightarrow (Y,\tau_d)$ is a quotient map. Therefore, from {\bf 2-d)}, the map $h:(Y,\tau_d)\rightarrow (\mathcal{D}_{f_i},\tau(\mathcal{D}_{f_i})),y\mapsto f_i^{-1}(y)$ must be a homeomorphism. Since $(Y,\tau_d)$ is self-similar, according to {\bf 2-e)}, the decomposition space $\mathcal{D}_{f_i}=\{f^{-1}_i(y)\subset X_i;y\in Y\}, \tau(\mathcal{D}_{f_i})=\{\mathcal{U}\subset \mathcal{D}_{f_i};\bigcup \mathcal{U}\in \tau_{X_i}\}$ of $(X_i,\tau_{X_i})$ is also self-similar. The space $X$ is filled with self-similar spaces in the sense that the family $\{f^{-1}_i(y);y\in Y,i=1,\dots,n\}$ of subsets of $X$ is a cover of $X$. The self-similar spaces $\mathcal{D}_{f_i}$ and $\mathcal{D}_{f_j}, i\not =j$ are disjoint in the sense that any point $f^{-1}_i(y)$ of $\mathcal{D}_{f_i}$ is not a point of $\mathcal{D}_{f_j}$ (Fig.1).

Finally, we note that in the above discussions we can replace the self-similar space with a compact substance in the materials science such as dendrite (A metric space is called a dendrite provided that it is a connected, locally connected, compact metric space which contains no simple closed curve as its subspace. A space which is homeomorphic to a dendrite is a dendrite. cf. {\bf 2-e)}.) 
[1,2]
 and then, we can obtain a family $\{X_1,\dots,X_n\}$ of subspace $X_i$ whose decomposition space is characterized as a dendrite.

\begin{figure}
\centering
\vspace{10mm}
\scalebox{.5}{\includegraphics[clip]{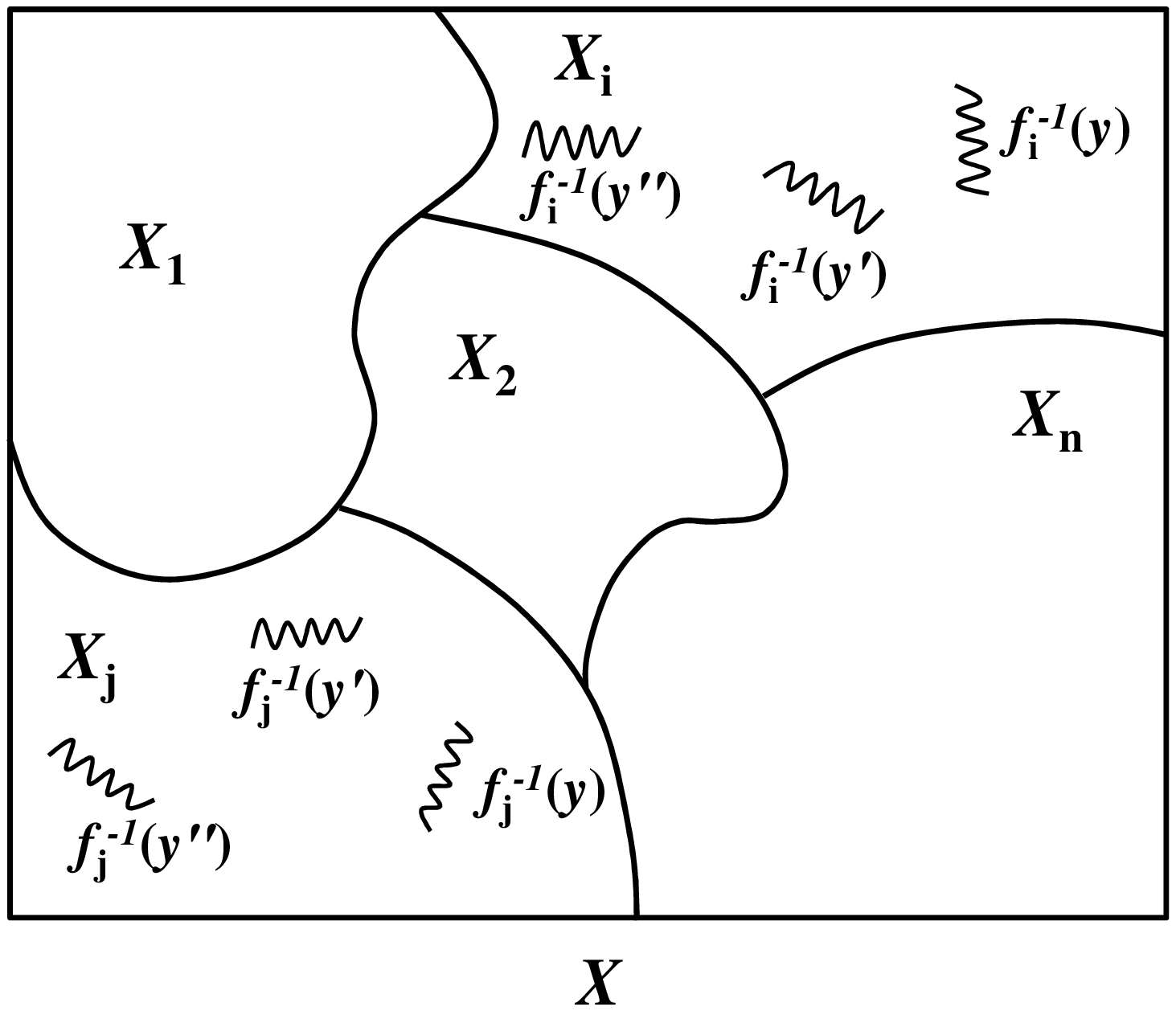}}\\
\caption{A schematic explanation of the mathematical procedure. $\{X_1,\dots,X_n\}$ is a partition of $X$. A decomposition space $\mathcal{D}_{f_i}=\{f^{-1}_i(y)\subset X_i;y\in Y\}$ of subspace $X_i$ of $X$ generated from a self-similar space $Y$ is self-similar. Each $f^{-1}_i(y)$ is a point of a self-similar space $\mathcal{D}_{f_i}$.}

\end{figure}
~~\\

\section{Conclusion}

Let a topological space $(X,\tau)$ be the $\Lambda$-product space $(\{0,1\}^\Lambda ,\tau_0^\Lambda )$ of $(\{0,1\},\tau_0)$. Here $Card \Lambda \succ \aleph$ (aleph) and $\tau_0=\{\{0,1\},\{0\},\{1\},\phi \}$. $(X,\tau)$ has an any $n$-partition $\{X_1,\dots,X_n\}$ where each subspace $(X_i,\tau_{X_i})$ is a non-empty clopen set of $(X,\tau)$ and has a self-similar decomposition space. Namely, $(X,\tau)$ is filled with the any $n$ mutually disjoint self-similar spaces each of which is a decomposition space of $(X_i,\tau_{X_i}), i\in\{1,\dots,n\}$. It is emphasized that the discussions concerning the self-similarity do not  depend on the metrizability of the initial space $(X,\tau)$.\\

\noindent
{\bf Acknowledgement}\\
The authors are grateful to Dr. Y. Yamashita, Dr. H. Ryo, Prof. Emeritus H.Fukaishi at Kagawa university, and Prof. H. Nagahama at Tohoku university for useful suggestions and encouragements.
 

\end{document}